\documentstyle{marcdek}

\makeindex

\def\sl{\em}

\def\Bbb{\bf}

\def\N{{\Bbb N}}
\def\C{{\Bbb C}}

\def\F{{\cal F}}

\def\invc{{c^{-1}}}

\def\listit#1,#2{#1_1$, $#1_2,\ldots,$\ $#1_{#2}}

\def\lnorm{\left\|}
\def\rnorm{\right\|}
\def\normo#1{\lnorm #1 \rnorm}
\def\widedot{\,\cdot\,}
\def\normdot{\normo{\widedot}}

\def\lmod{\left|}
\def\rmod{\right|}
\def\modo#1{\lmod #1 \rmod}

\def\Deltacond{$\Delta_2$-condition}
\def\phifunction{$\varphi$-function}
\def\Nfunction{{\it N}-function}
\def\conditionJ{condition~$(J)$}

\newcounter{rom}
\newenvironment{itemrom}{%
\begin{list}%
{\roman{rom})}%
{\usecounter{rom}%
\setlength\topsep{0in}\setlength\partopsep{0in}%
\setlength\itemsep{0in}\setlength\parsep{0in}}}%
{\end{list}}

\def\itemi{\item}
\def\itemii{\item}
\def\itemiii{\item}
\def\itemiv{\item}
\def\itemv{\item}

\newenvironment{defin}{\proclaim{DEFINITION}}{\endproclaim}
\newtheorem{thm}{THEOREM}[section]
\newtheorem{prop}[thm]{PROPOSITION}
\newtheorem{lemma}[thm]{LEMMA}

\begin{document}

\title{Boyd Indices of Orlicz--Lorentz Spaces}

\author{STEPHEN J.~MONTGOMERY-SMITH%
\thanks
{Research supported in part by N.S.F.\ Grants DMS 9001796 and DMS 9001357.}
\ \ \
Department of Mathematics, University of Missouri,
Columbia, Missouri 65211}

\maketitle

\thispagestyle{empty}
\pagestyle{empty}

\section*{ABSTRACT}

Orlicz--Lorentz spaces 
provide a common generalization of Orlicz spaces 
and
Lorentz spaces. In this paper, we investigate their Boyd indices. Bounds
on the Boyd indices in terms of the Matuszewska--Orlicz indices of the 
defining
functions are given. Also, we give an example to show that the Boyd indices
and Zippin indices of an Orlicz--Lorentz space need not be equal, answering a
question of Maligranda.  Finally, we show how the Boyd indices are related to
whether an Orlicz--Lorentz space is $p$-convex or $q$-concave.

\section{INTRODUCTION}

The Boyd indices of a rearrangement invariant space are of fundamental
importance. They were originally introduced by Boyd (1969) for the purpose
of showing certain interpolation results. Since then, they have played a 
major
role in the theory of rearrangement invariant spaces (see, for example,
Bennett and Sharpley (1988), Lindenstrauss and Tzafriri (1979) or 
Maligranda (1984)).

Orlicz--Lorentz spaces provide a common generalization of Orlicz spaces (see
Orlicz(1932) or Luxembourg (1955)) and 
Lorentz spaces\index{Lorentz space} (see Lorentz (1950) or 
Hunt (1966)), and have been
studied by many authors, including, for example, Maligranda (1984), Masty\l o
(1986) and Kami\'nska (1990a, 1990b, 1991). In particular,
Maligranda posed a question about the Boyd indices of these spaces. 

In this paper, we first give some
fairly elementary estimates for the Boyd indices of Orlicz--Lorentz 
spaces. Then we
give an example that show that
these estimates cannot be improved, thus answering 
Maligranda's question. Finally we show how knowledge of the Boyd indices
gives information about the $p$-convexity or $q$-concavity of the 
Orlicz--Lorentz
space.

\section{DEFINITIONS}

In discussing Orlicz--Lorentz spaces, it will be convenient to talk about 
them
in the more general framework of rearrangement invariant spaces. 
Unfortunately,
the definitions in the literature usually require that the spaces be
quasi-normed, which is not always the case with the Orlicz--Lorentz spaces. 
For
this reason we introduce the following definition of rearrangement invariant
spaces.

\begin{defin} If $(\Omega,\F,\mu)$\ is a measure space, we denote the 
measurable
functions, modulo functions equal to zero almost everywhere, by $L_0(\mu)$. 
We
say that a
{\em K\"othe functional}\ is a
function
$\normdot\colon L_0(\mu)\to[0,\infty]$\ satisfying
\begin{itemrom}
\itemi if $f\in L_0(\mu)$, then $\normo f=0\Leftrightarrow f=0$;
\itemii if $f\in L_0(\mu)$\ and $\alpha\in\C$, then $\normo{\alpha
f}=\modo\alpha
\normo f$;
\itemiii if $f,g\in L_0(\mu)$, then $\modo f\le\modo g \Rightarrow
\normo f\le
\normo g$;
\itemiv if $f_n,f\in L_0(\mu)$, then $\modo{f_n}\nearrow\modo
f \Rightarrow \normo{f_n}\to\normo f$;
\itemv if $f_n\in L_0(\mu)$, then $\normo{f_n}\to0
\Rightarrow f_n\to0$\ in the measure topology.
\end{itemrom}
\noindent
A {\em K\"othe space\/}\index{K\"othe space} 
is a pair $(X,\normdot)$, where $\normdot$\ is a
K\"othe functional, and $X=\{f\in L_0(\mu):\normo f <\infty\}$.
Usually, we will denote a space by a single letter, $X$, and denote its
functional by $\normdot_X$.
\end{defin}

\begin{defin} If $f:\Omega\to\C$\ is a measurable function, we define the 
{\em
non-increasing rearrangement\/}\index{non-increasing rearrangement} 
of $f$\ to be  
$$ f^*(x) = \sup\bigl\{\, t : \mu(\modo f \ge t) \ge x \,\bigr\} .$$ 
A {\em rearrangement invariant space\/}\index{rearrangement invariant space}
is a K\"othe
space such that if $f,g\in
L_0(\mu)$, and $f^* \le g^*$,
then $\normo f \le \normo g$.
\end{defin}

Now we define the Orlicz--Lorentz spaces. We refer the reader to 
Montgomery-Smith (1992) for a
motivation of the following definitions.

\begin{defin} A {\em \phifunction\/}\index{\phifunction} 
is a function $F:[0,\infty) \to
[0,\infty)$\ such that
\begin{itemrom}
\itemi $F(0) = 0$;
\itemii $\lim_{t\to\infty} F(t) = \infty$;
\itemiii $F$\ is strictly increasing;
\itemiv $F$\ is continuous;
\end{itemrom}

\noindent
We will say that a \phifunction\ $F$\ is {\em dilatory\/} if for some
$1<c_1,c_2<\infty$\ 
we have $F(c_1 t)\ge c_2 F(t)$\ for all $0\le t<\infty$. We
will say that $F$\ satisfies the {\em \Deltacond\/} if $F^{-1}$\ is dilatory.

If $F$\ is a \phifunction, we will define the function $\tilde F(t)$\ to
be $1/F(1/t)$\ if $t>0$, and $0$\ if $t=0$.

We say that two \phifunction s $F$\ and $G$\ are {\em
equivalent\/} (in symbols $F\asymp G$) if for some number $c<\infty$\ we have
that $F(\invc t) \le G(t) \le F(ct)$\ for all $0\le t<\infty$.

We will denote the \phifunction\ $F(t) = t^p$\ by $T^p$.
\end{defin}

\begin{defin} (See Orlicz (1932) or Luxembourg (1955).) 
If $(\Omega,\F,\mu)$\ is a measure
space, and $F$\ is a \phifunction, then we define the {\em Luxemburg
functional\/}\index{Luxembourg functional} of a measurable function $f$\ by
$$ \normo{f}_F = \inf\left\{\, c :
   \int_\Omega F\bigl(\modo{f(\omega)}/c\bigr) \,d\mu(\omega) 
   \le 1 \,\right\} 
   ,$$
The {\em Orlicz space\/}\index{Orlicz space} is the associated K\"othe space, 
and is denoted by
$L_F(\Omega,\F,\mu)$\ (or $L_F(\mu)$, $L_F(\Omega)$ or $L_F$\ for short).
\end{defin}

\begin{defin}
If $(\Omega,\F,\mu)$\ is a measure space, and
$F$\ and $G$\ are \phifunction s, then we define the {\em
Orlicz--Lorentz functional\/} of a measurable function $f$\ by 
$$ \normo{f}_{F,G} = \normo{f^*\circ\tilde F\circ\tilde G^{-1}}_G .$$
The {\em Orlicz--Lorentz space\/}\index{Orlicz--Lorentz space} 
is the associated K\"othe space, 
and is denoted
by $L_{F,G}(\Omega,\F,\mu)$\ (or $L_{F,G}(\mu)$, $L_{F,G}(\Omega)$ or 
$L_{F,G}$\
for short).

We will write $L_{F,p}$, $L_{p,G}$\ and $L_{p,q}$\ for $L_{F,T^p}$,
$L_{T^p,G}$\ and $L_{T^p,T^q}$\ respectively.
\end{defin}

It is an elementary matter to show that the Orlicz and Orlicz--Lorentz spaces
are rearrangement invariant spaces. We note that $\normdot_{F,F} = 
\normdot_F$,
and that $\normo{\chi_A}_{F,G} = \tilde F^{-1}\bigl(\mu(A)\bigr)$.

Now we define the various indices that we use throughout this paper.
Obviously, the most important of these are the Boyd indices. These were first
introduced in Boyd (1969). We will follow Maligranda (1984) for the names 
of the
other indices, but will modify the definitions so as to be consistent 
with the
notation used in Lindenstrauss and Tzafriri (1979). Thus other references 
to these indices will often
reverse the words `upper' and `lower', and use the reciprocals of the 
indices
used here. The Zippin
indices were introduced in Zippin (1971), and the Matuszewska--Orlicz 
indices in
Matuszewska and Orlicz (1960 and 1965). The Zippin indices are sometimes 
called {\em
fundamental indices}.

\begin{defin} For a rearrangement invariant space $X$, we let the 
{\em dilation
operators\/}\index{dilation operator} 
$d_a:X\to X$\ be $d_a f(x) = f(ax)$\  for $0<a<\infty$. 
We define
the {\em lower Boyd index\/}\index{Boyd index} to be
$$ p(X) =
   \sup \left\{\, p :
   \hbox{for some $c<\infty$\ we have $\normo{d_a}_{X\to X} \le c a^{-1/p}$\
         for $a<1$} \,\right\} .$$
We define the {\em upper Boyd index\/} to be
$$ q(X) =
   \inf \left\{\, q :
   \hbox{for some $c<\infty$\ we have $\normo{d_a}_{X\to X} \le c a^{-1/q}$\
         for $a>1$} \,\right\} .$$
We define the
{\em lower Zippin index\/}\index{Zippin index} to be
$$ p_z(X) =
   \sup \left\{\, p :
   \matrix{\hbox{for some $c<\infty$\ we have 
           $\normo{d_a\chi_A}_X \le c a^{-1/p}\normo{\chi_A}_X$} \cr
           \hbox{for all $a<1$\ and measurable $A$}\cr} 
   \,\right\} .$$ 
We define the {\em upper Zippin index\/} to be
$$ q_z(X) =
   \inf \left\{\, q :
   \matrix{\hbox{for some $c<\infty$\ we have 
           $\normo{d_a\chi_A}_X \le c a^{-1/q}\normo{\chi_A}_X$} \cr
           \hbox{for all $a>1$\ and measurable $A$}\cr} 
   \,\right\} .$$
\end{defin}

\begin{defin}
For a \phifunction\ $F$, we define the {\em lower Matuszewska--Orlicz 
index\/}\index{Matuszewska--Orlicz index}
to be
$$ p_m(F) =
   \sup \left\{\, p :
   \hbox{for some $c>0$\ we have $F(at) \ge c\, a^p F(t)$\
         for $0\le t<\infty$\ and $a>1$} \,\right\} .$$
We define the {\em upper Matuszewska--Orlicz index\/}
to be
$$ q_m(F) =
   \inf \left\{\, q :
   \hbox{for some $c<\infty$\ we have $F(at) \le c\, a^q F(t)$\
         for $0\le t<\infty$\ and $a>1$} \,\right\} .$$
\end{defin}

Thus, for example, 
$$ p(L_{p,q}) = q(L_{p,q}) = p_z(L_{p,q}) = q_z(L_{p,q})
   = p_m(T^p) = q_m(T^p) = p .$$
We also note the following elementary proposition about the 
Matuszewska--Orlicz
indices.

\begin{prop} Let $F$\ be a \phifunction. 
\begin{itemrom}
\itemi $F$\ is dilatory if and only if $p_m(F) > 0$.
\itemii $F$\ satisfies the \Deltacond\ if and only if $q_m(F) < \infty$.
\end{itemrom}
\end{prop}

It was conjectured, at one time, that the Boyd and Zippin indices coincide.
This is a natural conjecture in view of the fact that these indices do 
coincide
for almost all `natural' rearrangement spaces, for example, the Orlicz spaces
and the Lorentz spaces. However
Shimogaki (1970) gave an example of a rearrangement invariant Banach space
where these indices differ.

Maligranda (1984) posed a conjecture (Problem~6.1)
that would imply that the Boyd indices and Zippin indices coincide for the
Orlicz--Lorentz spaces. One of the main purposes of this paper is to show 
that
this is not the case.

In the sequel, we will always suppose that the measure space is $[0,\infty)$\
with the Lebsgue measure $\lambda$.

\section{ESTIMATES FOR THE BOYD INDICES OF THE ORLICZ LORENTZ SPACES}

The first results that we present give estimates for the Boyd indices. These
estimates are not very sophisticated. However, as we will show in Section~4,
they cannot be improved, at least in the form in which they are given. 
It would
be nice to give better estimates at some point in the future, which would 
make
use of more detailed structure information of the  defining functions of the
Orlicz--Lorentz space.

\begin{thm} Let $F$\ and $G$\ be \phifunction s. Then
\begin{itemrom}
\itemi
$ p_m(F) \ge p(L_{F,G}) \ge p_m(F\circ G^{-1}) p_m(G) \ge
p_m(F) p_m(G)/q_m(G)$; 
\itemii
$ q_m(F) \le q(L_{F,G}) \le q_m(F\circ G^{-1}) q_m(G) \le 
q_m(F) q_m(G)/p_m(G)$.
\end{itemrom}
\end{thm}

This will follow from the following propositions.

\begin{prop} Let $X$\ be a rearrangement invariant space, and let
$F$\ and $G$\ be \phifunction s.
\begin{itemrom}
\itemi $p(X) \le p_z(X)$\ and $q(X) \ge q_z(X)$.
\itemii $p(L_F) = p_z(L_F)$\ and $q(L_F) = q_z(L_F)$.
\itemiii $p_z(L_{F,G}) = p_m(F)$\ and $q_z(L_{F,G}) = q_m(F)$.
\end{itemrom}
\end{prop}

\begin{proof}
See Maligranda (1984) for part~(i), and see 
Lindenstrauss and Tzafriri (1979) for
part~(ii). Part~(iii) is clear.
\end{proof}

\begin{prop} Let $F_1$, $F_2$\ and $G$\ be \phifunction s. 
\begin{itemrom}
\itemi  $p(L_{F_1,G}) \ge p_m(F_1\circ F_2^{-1}) p(L_{F_2,G})$.
\itemii $q(L_{F_1,G}) \le q_m(F_1\circ F_2^{-1}) q(L_{F_2,G})$.
\end{itemrom}
\end{prop}

\begin{proof}
We will show (i). The proof of (ii) is similar. 

We note that if $p_1 < p_m(F_1\circ F_2^{-1})$, and if $p_2 < p(L_{F_2,G})$,
then there is a constant $c_1<\infty$\ such that for any $t\ge 0$\  and
$0<a<1$\ we have
$$ a \tilde F_1\circ\tilde F_2^{-1}(t)
   \le \tilde F_1\circ\tilde F_2^{-1}(c_1\,a^{1/p_1} t) ,$$
and there is a constant $c_2 < \infty$\ such that for any $f\in L_0$\ and
$0<b<1$\ we have
$$ \normo{d_{c_1 b}f}_{F_2,G} \le c_2 b^{-1/p_2} \normo f_{F_2,G} . $$
Therefore,
\begin{eqnarray*}
   \normo{d_a f}_{F_1,G}
   &=& \normo{ x\mapsto f^*\left( a\, 
   \tilde F_1\circ\tilde G^{-1}(x) \right)}_G \\
   &\le& \normo{ x\mapsto f^*\circ \tilde F_1\circ \tilde F_2^{-1}\left(
   c_1\,a^{1/p_1}
   \tilde F_2\circ\tilde G^{-1}(x) \right)}_G \\
   &\le& c_2 a^{-1/p_1 p_2} \normo f_{F_1,G} . 
\end{eqnarray*}
Therefore $p(L_{F_1,G}) \ge p_1 p_2$, and the result follows.
\end{proof}

\noindent
Proof of Theorem 3.1: The first inequality follows from Proposition~3.2. The
second inequality follows from Propositions~(3.2) and~(3.3). The third 
inequality
follows because
$$ p_m(F\circ G^{-1}) \ge p_m(F) p_m(G^{-1}) = p_m(F)/q_m(G) .$$
\endproof

\section{BOYD INDICES CAN DIFFER FROM ZIPPIN INDICES}

Now we show that Theorem~3.1 cannot be improved.
In so doing, we answer
Problem~6.1 posed by Maligranda (1984), by showing that the Boyd indices
and Zippin indices do not necessarily coincide for the Orlicz--Lorentz 
spaces.

\begin{thm} Given $0<p<q<\infty$, there is a 
\phifunction\ $G$\ such that $p_m(G) = p$, $q_m(G) = q$,
$ p(L_{1,G}) = p/q$, and $q(L_{1,G}) = q/p$.
\end{thm}

We also have the following interesting example, that shows that an
Orlicz--Lorentz space need not be quasi normed just because its defining
functions are dilatory.

\begin{thm} There is a dilatory \phifunction\ $G$\ such that
$L_{1,G}$\ is not a quasi-Banach space.
\end{thm}

At the heart of these results is the following lemma.

\begin{lemma} Suppose that $0<p,q<\infty$, $a>1$\ and $n_0, n_1 \in
\N$\ are such that
$$ (n_1-n_0) a^{-p}\left(1-a^{-(p+q)}\right) + a^{-2p-q} = 1 .$$
Suppose that $G$\ is a \phifunction\ such that for some $L,M>0$\ we have 
that 
\begin{eqnarray*}
   \tilde G(Ma^{2n}t) &=& La^{(p+q)n}t^p \\
   \tilde G(Ma^{2n+1}t) &=& La^{(p+q)n+p}t^q 
\end{eqnarray*}
for $1\le t\le a$\ and $n_0\le n\le n_1+1$. Then for all\/
$0\le\theta\le\inf\left\{q/p,1\right\}$, there are functions $f$\ and $g$\ 
such
that we have 
$$ \normo{d_{a^{-\theta}}f}_{1,G} 
   = a^{(q/p)\theta} \normo f_{1,G}
   \quad\hbox{and}\quad
   \normo{d_{a^{\theta}}g}_{1,G} 
   = a^{-(q/p)\theta} \normo g_{1,G} .$$
\end{lemma}

\begin{proof} We define the functions $f$\ and $g$\ by
\begin{eqnarray*}
   f(Mx) &=& \cases{M^{-1}a^{-2n_0-3} & if $0\le x<a^{2n_0}$ \cr
                        M^{-1}a^{-2n  -3} & if $a^{2n} \le x < a^{2n+2}$\    
                        and $n_0\le n\le n_1$ \cr
                        0 & if $a^{2n_1+2} \le x$ \cr} \\
   g(Mx) &=& \cases{M^{-1}a^{-2n_0-3-(p/q) \theta} 
                        & if $0\le x<a^{2n_0+\theta}$ \cr
                        M^{-1}a^{-2n  -3-(p/q) \theta} & if 
                        $a^{2n+\theta} \le x < a^{2n+2+\theta}$\    
                        and $n_0\le n\le n_1$ \cr
                        0 & if $a^{2n_1+2+\theta} \le x$ \cr} 
\end{eqnarray*}
so that $g = a^{-(p/q) \theta}d_{a^{-\theta}}f$. Then it is sufficient to 
show
that $\normo f_{1,G} = \normo g_{1,G} = 1$. We will only show that $\normo
g_{1,G} = 1$, as setting $\theta=0$ gives the other equality.

First, we note that if
$$ La^{(p+q)n+p\theta} \le x < La^{(p+q)(n+1)+p\theta} $$
then
$$ M a^{2n+\theta} \le \tilde G^{-1}(x) < M a^{2n+2+\theta} $$
and so
$$ g^*\circ\tilde G^{-1}(x) = M^{-1} a^{-2n-3-(p/q) \theta} $$
implying that
$$ G\circ g^*\circ\tilde G^{-1}(x) = 1/\tilde G(M a^{2n+3+(p/q) \theta}) 
   = 1/\left(La^{(p+q)n+2p+q+p\theta}\right) 
   = L^{-1}a^{-(p+q)n-2p-q-p\theta} .$$
Similarly, if $0\le x < La^{(p+q)n_0+p\theta}$, then 
$G\circ g^*\circ\tilde G^{-1}(x) = L^{-1}a^{-(p+q)n_0-2p-q-p\theta} $. Hence
\begin{eqnarray*}
   && \int_0^\infty G\circ g^*\circ\tilde G^{-1}(x) \,dx \\
   &=& \sum_{n=n_0}^{n_1}
      \int_{La^{(p+q)n+p\theta}}^{La^{(p+q)(n+1)+p\theta}} 
      G\circ g^*\circ\tilde G^{-1}(x) \,dx
      + \int_0^{La^{(p+q)n_0+p\theta}} 
      G\circ g^*\circ\tilde G^{-1}(x) \,dx \\
   &=& \sum_{n=n_0}^{n_1}
      \left( La^{(p+q)(n+1)+p\theta} - La^{(p+q)n+p\theta} \right) 
      L^{-1} a^{-n(p+q)-2p-q-p\theta} \\
   && \qquad
      + La^{(p+q)n_0+p\theta} L^{-1}a^{-n_0(p+q)-2p-q-p\theta} \\
   &=& (n_1-n_0) a^{-p}\left(1-a^{-(p+q)}\right) 
      + a^{-2p-q} \\
   &=& 1 ,
\end{eqnarray*}
as required.
\end{proof}

\noindent
Proof of Theorem 4.1:  Construct sequences of numbers $a_k$, $b_k$, $M_k$\ 
and
$N_k$\ ($k\ge 0$)\ such that $M_k$\ and $N_k$\ are integers, $a_k,b_k>0$,
\begin{eqnarray*}
   M_k a_k^{-p}\left(1-a_k^{-(p+q)}\right) + a_k^{-2p-q} &=& 1 ,\\
   N_k b_k^{-q}\left(1-b_k^{-(p+q)}\right) + b_k^{-p-2q} &=& 1 ,
\end{eqnarray*}
$a_k \to \infty$, and $b_k\to\infty$. Define sequences $A_k$\ and $B_k$\
inductively as follows: $A_0 = B_0 = 1$, $B_k = A_k a_k^{2M_k+2}$, and 
$A_{k+1}
= B_k b_k^{2N_k+2}$\ for $k\ge0$. Define $G$\ by
\begin{eqnarray*}
   G(1) &=& 1 ,\\
   G(A_k a_k^{2n}t) 
   &=& G(A_k) a_k^{(p+q)n} t^p \\
   G(A_k a_k^{2n+1}t) 
   &=& G(A_k) a_k^{(p+q)n+p} t^q \\
\noalign{\noindent for $0\le n\le M_k$\ and $1\le t\le a_k$,}
   G(B_k b_k^{2n}t) 
   &=& G(B_k) b_k^{(p+q)n} t^q \\
   G(B_k b_k^{2n+1}t) 
   &=& G(B_k) b_k^{(p+q)n+q} t^p \\
\noalign{\noindent for $0\le n\le N_k$\ and $1\le t\le b_k$, and}
   G(t) &=& \tilde G(t) 
\end{eqnarray*}
for $t<1$. Clearly $p_m(G) = p$\ and $q_m(G) = q$. From Lemma~4.3, we have 
that
$p(L_{1,G}) = p/q$\ and $q(L_{1,G}) = q/p $.
\endproof

\noindent
Proof of Theorem 4.2: Let $q=1$, and construct sequences of numbers $p_k$, 
$a_k$\
and $N_k$\ ($k\ge 0$)\ such that $N_k$\ is an integer, $a_k>0$,
$$ N_k a_k^{-p_k}\left(1-a_k^{-(p_k+q)}\right) + a_k^{-2p_k-q} = 1 , $$
$p_k\to\infty$, and $a_k^{q/p_k} \to \infty$. Define a sequence $A_k$\
inductively as follows: $A_0 = 1$, and $A_{k+1} = A_k a_k^{2N_k+2}$\ for
$k\ge0$. Define $G$\ by 
\begin{eqnarray*}
   G(1) &=& 1 ,\\
   G(A_k a_k^{2n}t) 
   &=& G(A_k) a_k^{(p_k+q)n} t^p_k \\
   G(A_k a_k^{2n+1}t) 
   &=& G(A_k) a_k^{(p_k+q)n+p_k} t^q \\
\noalign{\noindent for $0\le n\le M_k$\ and $1\le t\le a_k$, and}
   G(t) &=& \tilde G(t) 
\end{eqnarray*}
for $t>1$.
Then $p_m(G) = 1$. From Lemma~4.3, we have that $p(L_{1,G})
= 0$, and so by Theorem~5.3(ii) below, $L_{1,G}$\ cannot be a quasi-Banach
space. 
\endproof

\section{CONVEXITY AND CONCAVITY OF ORLICZ--LORENTZ SPACES}

An important property that one might like to know about K\"othe spaces is
whether it is $p$-convex or $q$-concave for some prescribed $p$\ or $q$. 
These
questions have already been settled for Orlicz spaces and Lorentz spaces. 

For
Lorentz spaces, it is almost immediate from their definition (Bennett and
Sharpley (1988) or Hunt (1966)) 
that $L_{p,q}$\ is $q$-convex if $p\ge q$, and $p$-concave if $p\le q$.
However, outside of these ranges, it is more difficult. In general, it is 
only
the case that $L_{p,q}$\ is $q\wedge(p-\epsilon)$-convex and
$p\vee(q+\epsilon)$-concave. These results are shown in many places, for 
example,
in Bennett and
Sharpley (1988) or Hunt (1966). For Orlicz--Lorentz spaces, the same 
methods of proof
work, and we present these results here.

First we define the notions of $p$-convexity and $q$-concavity. These notions
may also be found in, for example, Lindenstrauss and Tzafriri (1979).

\begin{defin} If $X$\ is a K\"othe space, we
say that $X$\ is {\em $p$-convex}\index{$p$-convex}, respectively
{\em $q$-concave}\index{$q$-concave}, if for some $C<\infty$\ we have
\begin{eqnarray*}
   \normo{\left(\sum_{i=1}^n \modo{f_i}^p\right)^{1/p}}_X
   &\le& C \, \left(\sum_{i=1}^n \normo{f_i}_X^p\right)^{1/p} ,\\
\noalign{\noindent respectively}
   \normo{\left(\sum_{i=1}^n \modo{f_i}^q\right)^{1/q}}_X
   &\ge& C^{-1} \left(\sum_{i=1}^n \normo{f_i}_X^q\right)^{1/q} ,
\end{eqnarray*}
for any $\listit f,n \in X$. 
\end{defin}

The most elementary result about $p$-concavity and $q$-convexity is the
following. This corresponds to the result that $L_{p,q}$\ is $q$-convex if
$p\ge q$, and $p$-concave if $p\le q$.

\begin{thm} Let $F$\ and $G$\ be \phifunction s.
\begin{itemrom}
\itemi If $G\circ T^{1/p}$\ is equivalent to a convex function and
$\tilde G\circ \tilde F^{-1}$\ is concave, then $L_{F,G}$\ is $p$-convex.
\itemii If $G\circ T^{1/q}$\ is equivalent to a concave function and
$\tilde G\circ \tilde F^{-1}$\ is convex, then $L_{F,G}$\ is $q$-concave.
\end{itemrom}
\end{thm}

\begin{proof} 
We will only prove (i), as the proof of (ii) is similar. We first use
the identity
$$ \normo f_{F\circ T^p,G\circ T^p} = \normo{\modo f^p}_{F,G}^{1/p} $$
to notice that without loss of generality we may take $p=1$. 

From Hardy, Littlewood and P\'olya (1952), Chapter~X, it follows that
$$ \normo{f}_{F,G} = \sup \normo{f\circ\sigma\circ\tilde F\circ\tilde 
   G^{-1}}_G ,$$
where the supremum is over all measure preserving maps
$\sigma:[0,\infty)\to[0,\infty)$. Since $G$\ is convex, it follows from
Krasnosel'ski\u\i\ and Ruticki\u\i\ (1961) that $\normdot_G$\ is $1$-convex. 
Now the result follows
easily.
\end{proof}

However, if we take the Boyd indices into account, we can also obtain the
following results.
These correspond to the result that says that $L_{p,q}$\ is
$q\wedge(p-\epsilon)$-convex and $p\vee(q+\epsilon)$-concave.

To state and prove these results, it is first necessary to recall
notation and results from Montgomery-Smith (1992).

\begin{defin} If $F$\ and $G$\ are \phifunction s, then say that $F$\ is {\em
equivalently less convex than\/} $G$\ (in symbols $F \prec G$) if $G\circ
F^{-1}$\ is equivalent to a convex function. We say that $F$\ is {\em
equivalently more convex than\/} $G$\ (in symbols $F \succ G$) if $G$\ is
equivalently less convex than $F$.

A \phifunction\ $F$\ is said to be an {\em \Nfunction\/}\index{\Nfunction}
if it is equivalent to a \phifunction\ $F_0$\ such that $F_0(t)/t$\ is
strictly increasing, $F_0(t)/t \to \infty$\ as $t\to \infty$, and $F_0(t)/t 
\to
0$\ as $t\to0$.

A \phifunction\ $F$\ is said to be {\em complementary\/}\index{complementary
function}
to a \phifunction\ $G$\ if for some $c<\infty$\ we have 
$$ \invc t \le F^{-1}(t) \cdot G^{-1}(t) \le ct 
   \qquad (0\le t<\infty) .$$
If $F$\ is an \Nfunction, we will let $F^*$\ denote a function complementary
to $F$.

An \Nfunction\ $H$\ is said to satisfy {\em \conditionJ\/}\index{\conditionJ} 
if
$$ \normo{1/\tilde H^*{}^{-1}}_{H^*} < \infty .$$
\end{defin}

To give some intuitive feeling for \Nfunction s that satisfy \conditionJ, 
we point out that these are functions
that equivalent to slowly rising convex functions, for example,
$$ F(t) = \cases{ t^{1+1/\log(1+t)} & if $t\ge1$\cr
                  t^{1-1/\log(1+1/t)} & if $t\le 1$.\cr}$$

\begin{thm} (Montgomery-Smith, 1992) Let $F$, $G_1$\ and $G_2$\ be
\phifunction s such that one of $G_1$\ and $G_2$\ is dilatory, and one of 
$G_1$\
or $G_2$\ satisfies the \Deltacond. Then the following are equivalent.
\begin{itemrom}
\itemi For some $c<\infty$, we have that 
$\normo f_{F,G_1} \le c\, \normo f_{F,G_2}$\ for all measurable $f$.
\itemii There is an \Nfunction\ $H$\ satisfying \conditionJ\ such that 
$G_1\circ
G_2^{-1} \succ H^{-1}$. 
\end{itemrom}
\end{thm}

Now, we are ready to state the main results of this section.

\begin{thm} Let $F$\ and $G$\ be \phifunction s, and $0<p<\infty$.
\begin{itemrom}
\itemi If the lower Boyd index
$p(L_{F,G}) > p$, and if $G\succ H^{-1}\circ T^p$\ for some \Nfunction\
satisfying \conditionJ, then $L_{F,G}$\ is $p$-convex.
\itemii If $L_{F,G}$\ is $p$-convex, then the lower Boyd index
$p(L_{F,G}) \ge p$, and $G\succ H^{-1}\circ T^p$\ for some \Nfunction\ 
satisfying
\conditionJ.
\end{itemrom}
\end{thm}

Note that in part~(i), it is not sufficient to take $p(L_{F,G}) = 1$. This is
shown by the example $L_{1,q}$\ for $1<q<\infty$, which is known to be not
1-convex (Hunt, 1966).

\begin{thm} Let $F$\ and $G$\ be \phifunction s such that $G$\ is
dilatory and $p(L_{F,G}) > 0$, and let $0<q<\infty$.
\begin{itemrom}
\itemi If the lower Boyd index
$q(L_{F,G}) < q$, and if $T^q \circ G^{-1} \succ H^{-1}$\ for some
\Nfunction\ satisfying \conditionJ, then $L_{F,G}$\ is $q$-concave.
\itemii If $L_{F,G}$\ is $q$-convex, then the lower Boyd index
$q(L_{F,G}) \le p$, and $T^q \circ G^{-1} \succ H^{-1}$\ for some \Nfunction\
satisfying \conditionJ.
\end{itemrom}
\end{thm}

\noindent Proof of Theorem 5.3: As in the beginning of the proof of 
Theorem~5.1, we may
suppose without loss of generality that $p=1$.

The proof of (i) uses fairly standard techniques 
(Bennett and Sharpley, 1988). First, by Theorem~5.2, 
we may assume that $G$\ is equivalent to 
a convex function. Next,
for any measurable function $f$, we define
$$ f^{**}(x) = {1\over x} \int_0^x f^*(\xi) \,d\xi
   = \int_0^1 d_a f^*(x) \,da .$$
Then we have the Hardy inequality holding, that is, for some $c<\infty$\ we
have that $\normo f_{F,G} \le \normo{f^{**}}_{F,G} \le c\,\normo f_{F,G}$. 
The
left hand inequality is obvious. For the right hand inequality, since
$p(L_{F,G}) > 1$, we know that for some $p>1$\ and some $c_1<\infty$\ we have
that $\normo{d_a}_{L_{F,G}\to L_{F,G}} \le c_1 a^{-1/p}$\ for all $a<1$. 
Hence
\begin{eqnarray*}
   \normo{f^{**}}_{F,G}
   &=& \normo{ \int_0^1 d_a f^* \,da }_{F,G} \\
   &=& \normo{ \int_0^1 d_a f^*\circ\tilde F\circ\tilde G^{-1} \,da }_G \\
   &\le& c_2 \int_0^1 \normo{ d_a f^*\circ\tilde F\circ\tilde 
         G^{-1} }_G \,da \\
\noalign{\noindent(as $G$\ is equivalent to a convex function)}
   &=& c_2 \int_0^1 \normo{ d_a f^* }_{F,G} \,da \\
   &\le& c_2 \int_0^1 c_1 a^{-1/p} \,da \, \normo{f^*}_{F,G} \\
   &\le& c_1 c_2 {p\over p-1} \normo f_{F,G} .
\end{eqnarray*}
But, the functional that takes $f$\ to $\normo{f^{**}}_{F,G}$\
is $1$-convex. This is because for any $x_0>0$, we have that
$$ f^{**}(x_0) = \sup_{\lambda(A) = x_0} \int_A \modo{f(x)}\,dx .$$
(See Hardy, Littlewood and P\'olya (1952), Chapter~X, or 
Lindenstrauss and Tzafriri (1979).)
Hence, $(f+g)^{**} \le f^{**} + g^{**}$. Also, by Krasnosel'ski\u\i\ and 
Ruticki\u\i\ (1961), 
it follows that $\normdot_G$\ is $1$-convex. Therefore,
\begin{eqnarray*}
   \normo{\sum_{i=1}^n \modo{f_i}}_{F,G}
   &\le& \normo{\left(\sum_{i=1}^n \modo{f_i}\right)^{**}}_{F,G} \\
   &\le& \normo{\sum_{i=1}^n f_i^{**}}_{F,G} \\
   &=& \normo{\sum_{i=1}^n f_i^{**}\circ \tilde F \circ\tilde G^{-1}}_G \\
   &\le& c_2 \sum_{i=1}^n \normo{f_i^{**}\circ \tilde F 
         \circ\tilde G^{-1}}_G \\
   &=& c_2 \sum_{i=1}^n \normo{f_i^{**}}_{F,G} \\
   &\le& c\,c_2 \sum_{i=1}^n \normo{f_i}_{F,G} ,
\end{eqnarray*}
as desired.

To show (ii), we note that if $a$\ is the reciprocal of an integer, then 
there
are functions $\listit g,{{a^{-1}}}$, with disjoint supports, and each 
with the
same distribution as $f$, so that $g_1+g_2+\ldots+g_{a^{-1}}$\ has the same
distribution as $d_a f$. Hence
$$ \normo{d_a f}_{F,G}
   \le c\, \left( \normo{g_1}_{F,G} + \normo{g_2}_{F,G}
   + \ldots + \normo{g_{a^{-1}}}_{F,G} \right)
   = c\, a^{-1} \normo f_{F,G} .$$
Hence $p(L_{F,G}) \ge 1$.

To show that $G\succ H^{-1}$\ for some \Nfunction\ satisfying
\conditionJ, we note the following inequalities.
\begin{eqnarray*}
   \normo f_{F,G} &=& \normo{\int_0^\infty 
                      \chi_{\modo f \ge t} \,dt}_{F,G} \\
            &\le& c\, \int_0^\infty \normo{\chi_{\modo f \ge t}}_{F,G} 
                  \,dt \\
            &=& c\, \int_0^\infty \tilde F^{-1}\bigl(
                    \mu\{\modo f \ge t\}\bigr)
                    \,dt \\
            &=& c\, \normo f_{F,1} .
\end{eqnarray*}
Now the result follows immediately from Theorem~5.2.
\endproof

\noindent
Proof of Theorem 5.4: As in the proof of Theorem~5.1, we may assume 
that $q=1$. 
To prove (i) we first note, by Theorem~5.2, we may assume that $G^{-1}$\ is
equivalent to  a convex function. Since $G$\ is dilatory, it follows that 
$G$\
is equivalent to a concave function (see Montgomery-Smith (1992), 
Lemma~5.5.2).

Next,
for any measurable function $f$, we define
$$ f_{**}(x) = f^*(x) + {1\over x} \int_x^\infty f^*(\xi) \,d\xi 
   = f^*(x) + \int_1^\infty d_a f^*(x) \,da .$$
Then, for some $c<\infty$\ we
have that $\normo f_{F,G} \le \normo{f_{**}}_{F,G} \le c\,\normo f_{F,G}$. 
The
left hand inequality is obvious. 

For the right hand inequality, we argue as
follows. Since $q(L_{F,G}) < 1$, we know that for some $q<1$\ and some
$c_1<\infty$\ we have that $\normo{d_a}_{L_{F,G}\to L_{F,G}} \le c_1 
a^{-1/q}$\
for all $a>1$. Since $G$\ is dilatory, it is easy to see that there is there
some $p>0$\ such that $G\circ T^{1/p}$\ is equivalent to a convex
function. Let $q<r<1$. Then there is a constant $c_2<\infty$, depending
upon $r$\ only, such that 
\begin{eqnarray*}
   f_{**}(x) 
   &\le& c_2 \left(\bigl(f^*(x)\bigr)^p + 
   {1\over x^{p/r}} \int_x^\infty \xi^{p/r-1} \bigl(f^*(\xi)\bigr)^p
   \,d\xi \right)^{1/p} \\
   &=& c_2 \left( \bigl(f^*(x)\bigr)^p + 
   \int_1^\infty a^{p/r-1} \bigl(d_a f^*(x)\bigr)^p
   \,da \right)^{1/p} .
\end{eqnarray*}
For if the right hand side is less than or equal to $1$, then it
is easily seen that
$$ f^*(\xi) \le 1 \wedge \left({x\over \xi -x}\right)^{1/r}
   \qquad (\xi > x) ,$$
and hence
$$ f_{**}(x) \le \int_1^\infty 1 \wedge (\theta-1)^{-1/r} \,d\theta .$$
Thus we have the following inequalities.
\begin{eqnarray*}
   \normo{f_{**}}_{F,G}
   &\le& c_2 \normo{ \left( (f^*)^p + \int_1^\infty a^{p/r-1} (d_a f^*)^p 
         \,da
         \right)^{1/p} }_{F,G} \\
   &=& c_2 \normo{ \left( (f^*\circ\tilde F\circ\tilde G^{-1})^p + 
          \int_1^\infty a^{p/r-1} 
          (d_a f^*\circ\tilde F\circ\tilde G^{-1})^p \,da \right)^{1/p}}_G \\
   &\le &c_3 \left( \normo{f^*\circ\tilde F\circ\tilde G^{-1}}_G^p +
            \int_1^\infty a^{p/r-1}
            \normo{ d_a f^*\circ\tilde F\circ\tilde G^{-1} }_G^p 
            \,da \right)^{1/p} \\
\noalign{\noindent(as $G\circ T^{1/p}$\ is equivalent to a convex function)}
   &=& c_3 \left( \normo{f^*}_{F,G}^p + 
          \int_1^\infty a^{p/r-1}
          \normo{ d_a f^*}_{F,G}^p 
          \,da \right)^{1/p} \\
   &\le& c_1 c_3 \left( 1 + \int_1^\infty a^{p/r-p/q-1} \,da \right)^{1/p}
        \normo f_{F,G} .
\end{eqnarray*}
But, the functional that takes $f$\ to $\normo{f_{**}}_{F,G}$\
is $1$-concave. This is because for any $x_0>0$, we have that
$$ f_{**}(x_0) = {1\over x_0} \int_0^\infty f(\xi) \,d\xi - f^{**}(x_0) .$$ 
Hence, $(f+g)_{**} \ge f_{**} + g_{**}$.
Also, by an argument similar to that given in M.A.~Krasnosel'ski\u\i\ and 
Ruticki\u\i\ (1961), 
it follows that $\normdot_G$\ is $1$-concave. Therefore,
\begin{eqnarray*}
   \sum_{i=1}^n \normo{f_i}_{F,G}
   &\le& \sum_{i=1}^n \normo{f_{i**}}_{F,G} \\
   &=& \sum_{i=1}^n \normo{f_{i**}\circ \tilde F \circ\tilde G^{-1}}_G \\
   &\le& c_3 \normo{\sum_{i=1}^n f_{i**}\circ \tilde F \circ\tilde 
         G^{-1}}_G \\
   &=& c_3 \normo{\sum_{i=1}^n f_{i**}}_{F,G} \\
   &\le& c_3 \normo{\left(\sum_{i=1}^n \modo{f_i}\right)^{**}}_{F,G} \\
   &\le& c\, c_1 \normo{\sum_{i=1}^n \modo{f_i}}_{F,G} ,
\end{eqnarray*}
as desired.

To show (ii), we note that if $a$\ is an integer, then there
are functions $\listit g,a$, with disjoint supports, and each with the
same distribution as $d_a f$, so that $g_1+g_2+\ldots+g_a$\ has the same
distribution as $f$. Hence
$$ \normo{f}_{F,G}
   \ge c^{-1} \left( \normo{g_1}_{F,G} + \normo{g_2}_{F,G}
   + \ldots + \normo{g_a}_{F,G} \right)
   = c^{-1} a \normo{d_a f}_{F,G} .$$
Hence $q(L_{F,G}) \le 1$.

To show that $G\prec H$\ for some \Nfunction\ satisfying
\conditionJ, we note the following inequalities.
\begin{eqnarray*}
   \normo f_{F,G} &=& \normo{\int_0^\infty \chi_{\modo f \ge t} 
                      \,dt}_{F,G} \\
            &\ge& c^{-1} \int_0^\infty \normo{\chi_{\modo f 
                  \ge t}}_{F,G} \,dt \\
            &=& c^{-1} \int_0^\infty \tilde F^{-1}\bigl(
                \mu\{\modo f \ge t\}\bigr)
                \,dt \\
            &=& c^{-1} \normo f_{F,1} .
\end{eqnarray*}
Now the result follows immediately from Theorem~5.2.
\endproof

\section{ADDITIONAL COMMENTS}

First, we remark that there is another definition of Orlicz--Lorentz spaces
given by Torchinsky (1976) (see also Raynaud (1990)). If $F$\ and
$G$\ are \phifunction s, then we define 
$$    \normo f_{F,G}^T
   = \inf \left\{\, c :
     \int_0^\infty G\bigl(\tilde F^{-1}(x) f^*(x)/c\bigr)
     \, {dx\over x} \le 1 \right\} ,$$
If $F$\ is dilatory and satisfy the \Deltacond, and if $G$\ is dilatory,
then it is very easy to calculate the Boyd indices of these spaces --- 
they are
precisely the same as their corresponding Matuszewska--Orlicz indices. This
follows from the fact that under these conditions, 
$ \normo{\chi_{[0,t]}}_{F,G}^T \approx \tilde F^{-1}(t) $\ 
(See Raynaud (1990) for more
details).

We also pose some questions.
\begin{itemrom}
\itemi What is the dual of an Orlicz--Lorentz space (when the space itself
is $1$-convex)? Is it another Orlicz--Lorentz space?
\itemii Is it possible to find more precise estimates for the Boyd indices of
Orlicz--Lorentz spaces?
\end{itemrom}
An approach to the last problem 
(at least for giving necessary and sufficient conditions for $p(L_{1,G}) =
q(L_{1,G}) = 1$\ is suggested in Montgomery-Smith (1991).

\section*{ACKNOWLEDGEMENTS}

This paper is an extension of work that I presented in my Ph.D.\ thesis
(1988). I would like to express my thanks to D.J.H.~Garling, my Ph.D.\
advisor, as well as the Science and Engineering Research Council who 
financed my
studies at that time.

I would also like to express gratitude to A.~Kami\'nska, W.~Koslowski and
N.J.~Kalton for their keen interest and useful conversations.

\newpage

\begin{theindex}

\item Lorentz space, 1
\item K\"othe space, 2
\item non-increasing rearrangement, 2
\item rearrangement invariant space, 2
\item \phifunction, 2
\item Luxembourg functional, 3
\item Orlicz space, 3
\item Orlicz--Lorentz space, 3
\item dilation operator, 3
\item Boyd index, 3
\item Zippin index, 3
\item Matuszewska--Orlicz index, 4
\item $p$-convex, 8
\item $q$-concave, 8
\item \Nfunction, 8
\item complementary function, 9
\item \conditionJ, 9

\end{theindex}
\end{document}